\documentclass[journal]{IEEEtran}
\usepackage[latin1]{inputenc}
\usepackage[cmex10]{amsmath}
\usepackage{amsfonts}
\usepackage{amssymb}
\usepackage{graphicx}
\usepackage{verbatim}
\usepackage{array}
\usepackage{multirow}
\usepackage{dcolumn}
\usepackage{color}
\usepackage[noadjust]{cite}
\usepackage{url}
\usepackage{balance}
\usepackage{booktabs}
\usepackage[mathscr]{euscript}
\DeclareGraphicsExtensions{.eps}
\usepackage{setspace}
\usepackage{graphicx,epstopdf}
\usepackage{comment}
\usepackage{enumerate}
\usepackage{flushend}

\usepackage{bbold}
\usepackage{fancyhdr,graphicx,amsmath,amssymb}
\usepackage{mathtools}

\usepackage{geometry}
\begin{document}
\title{Distribution System Operation {\color{black}Amidst} Wildfire-Prone Climate Conditions {\color{black}Under} Decision-Dependent Line Availability Uncertainty}

\author{ 
{\color{black}Alexandre Moreira, \IEEEmembership{Member,~IEEE,} Felipe~Pianc\'o,~\IEEEmembership{Student~Member,~IEEE,}} Bruno Fanzeres, \IEEEmembership{Member,~IEEE,} Alexandre Street, \IEEEmembership{Senior Member,~IEEE,}  Ruiwei Jiang, Chaoyue Zhao, and Miguel Heleno, \IEEEmembership{Member,~IEEE}

}

\maketitle
\begin{abstract}
    Wildfires can severely damage electricity grids leading to long periods of power interruption. Climate change will exacerbate this threat by increasing the frequency of dry climate conditions. Under these climate conditions, human-related actions that initiate wildfires should be avoided, including those induced by power systems operation. In this paper, we propose a novel optimization model that is capable of determining appropriate network topology changes (via switching actions) to alleviate the levels of power flows through vulnerable parts of the grid so as to decrease the probability of wildfire ignition. Within this framework, the proposed model captures the relationship between failure probabilities and line-flow decisions by explicitly considering the former as a function of the latter. The resulting formulation is a two-stage model with endogenous decision-dependent probabilities, where the first stage determines the optimal switching actions and the second stage evaluates the worst-case expected operation cost. We propose an exact iterative method to deal with this intricate problem and the methodology is illustrated with a 54-bus and a 138-bus distribution system.
\end{abstract}


\vspace{-0.4cm}
\section*{Nomenclature}\label{Nomenclature}
\vspace{-0.4cm}
\subsection*{Sets}
\begin{description} [\IEEEsetlabelwidth{5000000}\IEEEusemathlabelsep]
	
    \item[$\mathcal{L}$] Set of indexes of line segments.
    \item[$\mathcal{L}^{sw}$] Set of indexes of switchable line segments.
    \item[$\mathcal{K}^{forbid}$] Set of indexes of forbidden switching patterns.
    \item[$\mathcal{N}$] Set of indexes of buses.
    \item[$\mathcal{N}^{subs}$] Set of indexes of buses with substation.
    
\end{description}

\vspace{-0.4cm}
\subsection*{Parameters}
\begin{description} [\IEEEsetlabelwidth{5000000}\IEEEusemathlabelsep]
	
    \item[$\beta_l$]{Sensitivity of failure probability to the scheduled active power flow of line $l \in \mathcal{L}$.}
    \item[$\gamma_l$]{Estimated upper bound for the nominal probability of failure associated with line $l \in \mathcal{L}$.}
    \item[$C^{ll}$]{Cost of loss of load.}
    \item[$C^{sw}_l$]{Cost of switching line $l \in \mathcal{L}^{sw}$.}
    \item[$C^{tr}_b$]{Cost of active power from main transmission grid for bus $b \in \mathcal{N}^{subs}$.}
    \item[$D^p_b$]{Active power demand at bus $b \in \mathcal{N}$.}
    \item[$E_l$]{Number of digits for binary expansion used in Master problem linearization for line $l \in \mathcal{L}$.}
    \item[$\overline{F}_l$]{Maximum power flow at line $l \in \mathcal{L}$.}
    \item[$\overline{P}^{tr}_b$]{Maximum active power injection at bus $b \in \mathcal{N}^{subs}$.}
    \item[$PF_b$]{Power factor at bus $b \in \mathcal{N}$.}
    \item[$\overline{Q}^{tr}_b$]{Maximum reactive power at bus $b \in \mathcal{N}^{subs}$.}
    \item[$\underline{Q}^{tr}_b$]{Minimum reactive power at bus $b \in \mathcal{N}^{subs}$.}
    \item[$R_l$]{Resistance of line $l \in \mathcal{L}$.}
    \item[$s$]{Step for binary expansion used in Master problem linearization.}
    \item[$\underline{V}_b$]{Voltage lower bound at bus $b \in \mathcal{N}$.}
    \item[$\overline{V}_b$]{Voltage upper bound at bus $b \in \mathcal{N}$.}
    \item[$V^{ref}$]{Voltage reference.}
    \item[$X_l$]{Reactance of line $l \in \mathcal{L}$.}
    \item[$z_l^{sw,0}$]{Initial switching status of line $l \in \mathcal{L}^{sw}$.}
\end{description}

\vspace{-0.6cm}
\subsection*{Decision variables}
\begin{description} [\IEEEsetlabelwidth{5000000}\IEEEusemathlabelsep]

    \item[$\alpha$]{Worst expected value of lower-level problem.}
    \item[$\Delta D^{p-}_b$]{Amount of active power loss at bus $b \in \mathcal{N}$.}
    \item[$\Delta D^{p+}_b$]{Amount of active power surplus at bus $b \in \mathcal{N}$.}
    \item[$\Delta D^{q-}_b$]{Amount of reactive power loss at bus $b \in \mathcal{N}$.}
    \item[$\Delta D^{q+}_b$]{Amount of reactive power surplus at bus $b \in \mathcal{N}$.}
    
    \item[$f^{p}_l$]{Active power flow at line $l \in \mathcal{L}$.}
    \item[$f^{q}_l$]{Reactive power flow at line $l \in \mathcal{L}$.}
    \item[$p^{tr}_b$]{Amount of active power injected at bus $b \in \mathcal{N}^{subs}$.}
    \item[$q^{tr}_b$]{Amount of reactive power at bus $b \in \mathcal{N}^{subs}$.}
    \item[$v^{\dagger}_b$]{Squared voltage at bus $b \in \mathcal{N}$.}
    \item[$y^{sw}_l$]{Binary decision variable indicating a switching action of line $l \in \mathcal{L}^{sw}$ (1 if switched, 0 otherwise).}
    \item[$z^{sw}_l$]{Binary decision variable of switching status of line $l \in \mathcal{L}^{sw}$ (1 if switched on, 0 otherwise).}
    
\end{description}

%
\vspace{-0.4cm}
\section{Introduction}
\vspace{-0.2cm}
\IEEEPARstart{W}{ildfire} events are a real threat to power systems operations at both transmission and distribution levels. The damage caused by these events might cost a significantly large amount of irrecoverable capital to society (e.g., the estimate of more than \$700 million in damage to transmission and distribution systems over 2000-2016 \cite{WildFireImpact2018}) and be irreparable in cases when human lives are involved. Over the past two decades, California, for instance, has experienced a large raise in the frequency of small wildfires, while the total burned area from large ones has also substantially increased \cite{Li2019}. In this context, human-induced activities have been placed at a top rank among the main roots of wildfire ignition, with power system operations responsible for some of them, as, for instance, when eventual sparks due to power flow through overhead lines aligned with dry weather conditions and strong wind speed levels cause this natural disaster \cite{Russell2012}. In extreme cases, this has been addressed by the electric sector with public safety power shut-offs (PSPS), which results in significant load sheddings and economic impacts \cite{Chiu2022}. As a consequence, novel operative policies are of significant importance in order to establish efficient power system operations amidst wildfire-prone climate conditions, assuring thus high levels of sustainability and system resilience \cite{Davoudi2021_ReclosingDistSyst_WildFire, Nazaripouya2020_PowerGridRes_WildFire}.

Due to this critical prospect, various research efforts have been dedicated to addressing resilience in power systems under potential natural disasters and human-made attacks. At the transmission level, for example, the work developed in \cite{Romero2013} proposes a two-stage stochastic Mixed-Integer NonLinear Programming (MINLP) model to define investment strategies to improve resilience, considering a range of earthquake events and the methodology developed by \cite{Lagos2020} combines optimization and simulation techniques to determine a portfolio of investments to improve grid resilience while also considering the potential occurrence of earthquakes. At the distribution level, on the other hand, the work reported in \cite{Nazemi2020} presents a storage sitting and sizing model to increase resilience while facing seismic hazards and, in \cite{Lin2018}, the authors designed a three-level system of optimization models to identify line hardening solutions to protect the distribution grid against intentional or unintentional attacks. Particularly regarding wildfires, an increasing deal of attention has been emerging in technical literature. Notably, in \cite{Rhodes2021}, the authors propose a methodology to alleviate wildfire risks by optimizing the selection of components in the grid to be de-energized in a power shut-off scheme. In addition, in \cite{Trakas2018}, a stochastic programming model that aims at increasing the resiliency of a distribution system exposed to an approaching wildfire is devised under exogenous uncertainties such as solar radiation, wind speed, and wind direction. Notwithstanding the relevance of recent technical literature, none of them has taken into account the direct impact of the power flow dispatch on the likelihood of line failures in a decision-dependent uncertainty framework.

From a modeling perspective, it is important to emphasize that uncertainties in power system operations are typically exogenously induced into the decision-making process. In this framework, uncertainty sources are solely associated with external factors and are not endogenously affected by operational actions. However, in many realistic cases, such as under wildfire-prone climate conditions, the operation of electric grids is also associated with the origin of fire ignitions, which significantly increase line failure probabilities. Due to this double role of power grids, the nature of the uncertainty is thus more complex to characterize (dependent not only on meteorological conditions -- exogenous factors, but also on the grid operation decisions -- endogenous factors), challenging the standard exogenously-induced approach. Therefore, in order to design resilience-oriented operational strategies in high fire-threat areas, utility operators must be aware of the impact of their operational decisions on the likelihood of wildfire initiation and reduction in reliability levels\footnote{We refer to \cite{Muhs2020} and the references therein for a wider discussion on the impact of distribution system operations in the probabilistic characterization of wildfire ignition}, i.e., the endogenous nature of the uncertainty.

Methodologically, we can divide decision-dependent uncertainty characterizations into two major types. The first one involves problems where the decision-making process filters the potential paths of uncertainty realization \cite{Goel2004}. This filtering implies that a given decision might not only rule out possible scenarios from occurring but also open the possibility for a specific subset of future events to occur. The second type is associated with problems where the decisions directly impact the whole probabilistic characterization of the uncertainty factors \cite{Zhan2017_ProbabilityDDU}.

In this paper, we leverage the second modeling type to propose a new methodology for distribution system operations capable of endogenously taking into account the impact of power flows on failure probabilities in the context of a potential wildfire event. We design a decision-dependent uncertainty framework where the line failure probabilities are a function (dependent) of the power flow levels. In this framework, we consider that during adverse climate conditions (dry weather and reasonably strong wind), switching actions can be made to reduce power flows in vulnerable areas of the grid, therefore decreasing the probability of wildfire ignition and consequent line failures, while seeking to maintain load supply. Thus, the proposed methodology allows distribution system operators to perform efficient switching actions to improve the system reliability accounting for decision-dependent line availability uncertainty. Structurally, the proposed methodology falls into the class of a two-stage, distributionally robust optimization problem with decision-dependent uncertainty \cite{Luo2020_D3RO}. In the first stage, our model decides the network topology (switching lines) and power imports from the main grid with main goal of minimizing the operational cost in the pre-contingency state plus the worst-case expected cost of operating the system under post-contingency states considering probabilities adjusted according to the pre-contingency network topology and line power flows. Then, in the second stage, the power flow and energy not served are evaluated for each post-contingency state. To summarize, the contributions of this paper are twofold:

\begin{enumerate}
	
	\item To formulate the distribution grid operation under adverse climate conditions as a two-stage distributionally robust optimization problem where the probabilities of line failure are co-dependent on the weather conditions (exogenous) and system power flows (endogenous). In the first stage, the system operator decides switching actions (therefore determining grid topology) and power imports from the main grid aiming at co-optimizing the pre-contingency and the worst-case expected post-contingency operations, formulated as the second stage.
	
	\item To devise an effective decomposition-based solution methodology capable of solving the proposed optimization problem. The approach is able to circumvent the computational difficulties posed by the multi-level (non-convex) structure intrinsic to the decision-dependent uncertainty modeling frameworks.
	
\end{enumerate}
\vspace{-0.4cm}
\section{Optimal Distribution System Operation with Decision-Dependent Uncertainty} \label{MathematicalFormulation}

The main objective of this work is to propose a methodology to determine the least-cost operation of a distribution system taking into account the impact of operative decisions in the probabilistic characterization of the line availability. We assume that the operator can perform switching actions in a set of line segments in the distribution system with the objective of minimizing the worst-case expected operation cost considering a decision-dependent uncertainty in line availability. In \eqref{eq_model4_upper_primal00}--\eqref{eq_model4_upper_primal34}, the proposed distribution system operation model is formulated.

\begin{align}
    & \underset{{\substack{\Delta D_b^{p-}, \Delta D_b^{p+}, \Delta D_b^{q-}, \Delta D_b^{q+}, \\
                           f^{p}_l,
                           f^{q}_l,
                           p^{{\color{black} tr }}_b, 
                           q^{{\color{black} tr }}_b,
                           v^{\dagger}_b,
                           y^{sw}_l,
                           z^{sw}_l}}}{\text{Minimize}} 
    ~ \hspace{0.50cm} \sum_{b\in \mathcal{N}^{subs}}{\Big( C^{tr}_b} p^{{\color{black}tr}}_b \Big) \notag \\
    & \hspace{0.50cm} + \sum_{b \in \mathcal{N}} C^{\color{black}{ll}} \Bigl(  \Delta D_b^{p+} + \Delta D_b^{p-} + \Delta D_b^{q+} + \Delta D_b^{q-} \Bigl) \notag \\
    & \hspace{0.50cm} + \sum_{l \in {\mathcal{L}}^{sw}} C^{sw}_l y^{sw}_l + \sup_{{\mathcal Q} \in {\mathcal P}(\boldsymbol{f}^{p}, \boldsymbol{\beta})}\mathbb{E}_{\mathcal Q} \Bigl[ H\bigl(\boldsymbol{z}^{sw},\boldsymbol{a}^L\bigr) \Bigr]  \label{eq_model4_upper_primal00} \\ \notag
    & ~ \text{subject to:} \notag \\ 
    & p^{{\color{black} tr }}_b + \sum_{l \in {\mathcal L}|to(l)=b} f_l^{p} - \sum_{l \in {\mathcal L}|fr(l)=b} f_l^{p} - D^{p}_b \notag \\
    & \hspace{74pt} - {\Delta}D^{p+}_b + {\Delta}D^{p-}_b = 0; 
    \forall b \in \mathcal{N}^{subs} \label{eq_model4_upper_primal01} \\
    & q^{{\color{black} tr }}_b + \sum_{l \in {\mathcal L}|to(l)=b} f^{q}_l - \sum_{l \in {\mathcal L}|fr(l)=b} f^{q}_l  \notag \\
    & \hspace{10pt} - \tan(\arccos(PF_b)) D^{p}_b - {\Delta}D^{q+}_b + {\Delta}D^{q-}_b = 0;\notag \\
    & \hspace{177pt} 
    \forall b \in \mathcal{N}^{subs} \label{eq_model4_upper_primal02} \\
    & \sum_{l \in {\mathcal L}|to(l)=b} f_l^{p} - \sum_{l \in {\mathcal L}|fr(l)=b} f_l^{p} - D^{p}_b - {\Delta}D^{p+}_b  \notag \\
    & \hspace{95pt} + {\Delta}D^{p-}_b = 0; 
    \forall b \in \mathcal{N} \setminus \mathcal{N}^{subs} \label{eq_model4_upper_primal03} \\
    & \sum_{l \in {\mathcal L}|to(l)=b} f^{q}_l - \sum_{l \in {\mathcal L}|fr(l)=b} f^{q}_l - \tan(\arccos(PF_b)) D^{p}_b \notag \\
    & \hspace{55pt} - {\Delta}D^{q+}_b + {\Delta}D^{q-}_b = 0; 
    \forall b \in \mathcal{N} \setminus \mathcal{N}^{subs} \label{eq_model4_upper_primal04} \\
    & - v^{\dagger}_{fr(l)} + v^{\dagger}_{to(l)} + 2(R_l f^{p}_l + X_l f^{q}_l) - (1 - z^{sw}_l)M \leq 0; \notag \\
    & \hspace{185pt}  
    \forall l \in \mathcal{L}^{sw} \label{eq_model4_upper_primal05} \\
    & v^{\dagger}_{fr(l)} - v^{\dagger}_{to(l)} - 2(R_l f^{p}_l + X_l f^{q}_l) - (1 - z^{sw}_l)M \leq 0; \notag \\
    & \hspace{185pt} 
    \forall l \in \mathcal{L}^{sw} \label{eq_model4_upper_primal06} \\
    & v^{\dagger}_{fr(l)} - v^{\dagger}_{to(l)} - 2(R_l f^{p}_l + X_l f^{q}_l) = 0; 
    \forall l \in \mathcal{L} \setminus \mathcal{L}^{sw} \label{eq_model4_upper_primal07} \\
    & \underline{V}_b^2 \leq v^{\dagger}_{b} \leq \overline{V}_b^2; 
    \forall b \in \mathcal{N} \label{eq_model4_upper_primal08} \\
    %
    %
    & v^{\dagger}_b = V^{ref^2} ; 
    \forall b \in \mathcal{N}^{subs} \label{eq_model4_upper_primal10} \\ 
    & - z^{sw}_l \overline{F}_l \leq f^{p}_l \leq z^{sw}_l \overline{F}_l; 
    \forall l \in {\mathcal L}^{sw} \label{eq_model4_upper_primal11} \\   
    %
    %
    & - z^{sw}_l \overline{F}_l \leq f^{q}_l \leq z^{sw}_l \overline{F}_l; 
    \forall l \in {\mathcal L}^{sw} \label{eq_model4_upper_primal13} \\
    %
    %
    & - \overline{F}_l \leq f^{p}_l \leq \overline{F}_l; 
    \forall l \in {\mathcal L} \label{eq_model4_upper_primal15} \\
    %
    %
    & - \overline{F}_l \leq f^{q}_l \leq \overline{F}_l; 
    \forall l \in {\mathcal L} \label{eq_model4_upper_primal17} \\
    %
    %
    & f^{q}_l - {\color{black}\cot} \Biggl ( \Biggl ( \frac{1}{2} - e  \Biggr ) \frac{\pi}{4} \Biggr ) \Biggl ( f^{p}_l - {\color{black}\cos} \Biggl ( e \frac{\pi}{4} \Biggr ) \overline{F}_l \Biggr ) \notag \\
    & \hspace{41pt} - {\color{black}\sin} \Biggl ( e \frac{\pi}{4} \Biggr )\overline{F}_l \leq 0; 
    \forall l \in {\mathcal L}, e \in \{1,\dots,4\} \label{eq_model4_upper_primal19} \\
    & - f^{q}_l - {\color{black}\cot} \Biggl ( \Biggl ( \frac{1}{2} - e  \Biggr ) \frac{\pi}{4} \Biggr ) \Biggl ( f^{p}_l - {\color{black}\cos} \Biggl ( e \frac{\pi}{4} \Biggr ) \overline{F}_l \Biggr ) \notag \\
    & \hspace{41pt} - {\color{black}\sin} \Biggl ( e \frac{\pi}{4} \Biggr ) \overline{F}_l \leq 0; 
    \forall l \in {\mathcal L}, e \in \{ 1,\dots,4 \} \label{eq_model4_upper_primal20} \\
    & 0 \leq p^{{\color{black} tr }}_b \leq \overline{P}_b^{tr}; 
    \forall b \in \mathcal{N}^{subs} \label{eq_model4_upper_primal21} \\
    %
    %
    & \underline{Q}_b^{tr} \leq q^{{\color{black} tr }}_b \leq \overline{Q}_b^{tr}; 
    \forall b \in \mathcal{N}^{subs} \label{eq_model4_upper_primal23} \\
    %
    %
    & {\Delta}D^{p+}_{b}, {\Delta}D^{p-}_{b}, {\Delta}D^{q+}_{b} , {\Delta}D^{q-}_{b} \geq 0; \forall b \in \mathcal{N}\label{eq_model4_upper_primal24} \\
    %
    %
    %
    %
    %
    & \Delta D_b^{p^{-}} \leq D^p_b; 
    \forall b \in \mathcal{N} \label{eq_model4_upper_primal29} \\
    & \Delta D_b^{q^{-}} \leq \tan(\arccos(PF_b))(D^p_b); 
    \forall b \in \mathcal{N} \label{eq_model4_upper_primal30} \\
    & y^{sw}_l \geq  z^{sw}_l - z^{sw,0}_l; 
    \forall l \in {\mathcal L}^{sw} \label{eq_model4_upper_primal31} \\
    & y^{sw}_l \geq z^{sw,0}_l - z^{sw}_l; 
    \forall l \in {\mathcal L}^{sw} \label{eq_model4_upper_primal32} \\ 
    & {\color{black}\sum_{l \in {\mathcal L}^{forbid}_k} z^{sw}_l \leq \Bigl |{\mathcal L}^{forbid}_k \Bigr | - 1; 
    \forall k \in {\mathcal K}^{forbid} } \label{eq_model4_upper_primal33} \\ 
    & z^{sw}_l \in \{0,1\}; 
    \forall l \in {\mathcal L}^{sw}\label{eq_model4_upper_primal34}
\end{align}


\noindent where sets ${\cal L}$, ${\cal L}^{sw}$, ${\cal K}^{forbid}$, ${\cal N}$, and ${\cal N}^{subs}$ contain indices of all line segments, line segments that can be switched on/off, line segments that cannot be simultaneously switched on (due to radiality constraints), all buses of the distribution system, and buses with substations, respectively. In addition, parameters $C^{tr}_b$, $C^{\color{black}{ll}}$, $C^{sw}_l$, $z_l^{sw,0}$, $D^p_b$, $PF_{b}$, $V^{ref}$, $R_l$, $X_l$, ${\color{black}\underline{V}_b}$, ${\color{black}\overline{V}_b}$, {\color{black}{$\overline{F}_l$}}, $\overline{P}^{{\color{black} tr }}_b$, $\overline{Q}^{{\color{black} tr }}_b$, ${\color{black}\underline{Q}^{{\color{black} tr }}_b}$ represent cost of purchasing active power from the main transmission grid, cost of loss of load, cost of switching, initial switching status of switchable line segments (equal to 1 if switched on, 0 otherwise), active power demand, power factor, voltage reference, resistance, reactance, voltage lower bound, voltage upper bound, {\color{black}maximum power flow in each line segment}, maximum active power injection at the substations, maximum reactive power injection at the substations, {\color{black}and minimum reactive power injection at the substations}, respectively. Moreover, decision variables $p^{{\color{black} tr }}_{b}$, $q^{{\color{black} tr }}_{b}$, $v^{\dagger}_{b}$, $f^p_{l}$, $f^q_{l}$, $y^{sw}_{l}$, $z^{sw}_{l}$, ${\Delta}D^{p+}_{b}$, ${\Delta}D^{p-}_{b}$, ${\Delta}D^{q+}_{b}$, ${\Delta}D^{q-}_{b}$ represent active power injected at the substations, reactive power injected at the substations, squared voltage, active power flow, reactive power flow, an indication of a switching action (equal to 1 if a switching action is scheduled, 0 otherwise), switching status, active power surplus, active power loss, reactive power surplus, reactive power loss.

Problem \eqref{eq_model4_upper_primal00}--\eqref{eq_model4_upper_primal34} is a two-stage, mixed-integer, distributionally robust optimization problem with decision-dependent uncertainty (ambiguity set). The objective function \eqref{eq_model4_upper_primal00} aims at minimizing a combination of active power injection purchases at the nodes with substations, loss of load costs, switching action, as well as the decision-dependent expected second-stage operational cost. More specifically, the latter is represented by $H(\mathbf{z}^{sw},\boldsymbol{a}^{L})$, a function of the first stage switching decision $(\mathbf{z}^{sw})$ and the random vector $\boldsymbol{a}^L$ associated with the availability of line segments of the feeder. To do so, note that in \eqref{eq_model4_upper_primal00}, we formulate the ambiguity set ${\cal P}$ (that accounts for the collection of credible probability distributions of line availability uncertainty - this set will be better defined in Subsection \ref{sec::AmbiguitySet}) as a function of the scheduled power flow $(\boldsymbol{f}^{p})$ and the (contextual) factors $(\boldsymbol{\beta})$ (also better defined in Subsection \ref{sec::AmbiguitySet}) to characterize endogenous and exogenous influence to uncertainty in line failures, respectively.


Active and reactive power balance are modeled through constraints \eqref{eq_model4_upper_primal01} and \eqref{eq_model4_upper_primal02} for substations and via constraints \eqref{eq_model4_upper_primal03} and \eqref{eq_model4_upper_primal04} for the remaining buses. Constraints \eqref{eq_model4_upper_primal05} and \eqref{eq_model4_upper_primal06} model voltage difference between sending and receiving ends of switchable line segments, with $M$ denoting a large number to relax these constraints when line $l \in {\cal L}^{sw}$ is switched off. Analogously, constraints \eqref{eq_model4_upper_primal07} represent voltage drop for non-switchable line segments. Constraints \eqref{eq_model4_upper_primal08} enforce voltage limits. Constraints \eqref{eq_model4_upper_primal10} set the voltage at the substations equal to the voltage reference. Active power flows limits are imposed by constraints \eqref{eq_model4_upper_primal11} for switchable line segments and by \eqref{eq_model4_upper_primal15} for the remaining ones. Likewise, constraints \eqref{eq_model4_upper_primal13} and \eqref{eq_model4_upper_primal17} impose limits to reactive power flows, which are also limited according to current active power flows by constraints \eqref{eq_model4_upper_primal19} and \eqref{eq_model4_upper_primal20} similarly to the linearized AC power flow presented in \cite{Mashayekh2018}. Constraints \eqref{eq_model4_upper_primal21} and \eqref{eq_model4_upper_primal23} enforce limits on active and reactive power injections at the substations, respectively. Constraints \eqref{eq_model4_upper_primal24} enforce non-negativity to power surplus and load shedding variables while constraints \eqref{eq_model4_upper_primal29} and \eqref{eq_model4_upper_primal30} impose upper limits on load shedding variables. Constraints \eqref{eq_model4_upper_primal31} and \eqref{eq_model4_upper_primal32} model the behavior of variable $y^{sw}_{l}$, which assumes value equal to 1 if the determined switching status $z^{sw}_{l}$ of line segment $l \in {\cal L}^{sw}$ is different from its initial switching status $z^{sw,0}_{l}$. Constraints \eqref{eq_model4_upper_primal33} model the forbidden switching patterns with ${\cal L}^{forbid}_{k}$ indicating the lines segments that cannot be simultaneously switched on for each $k \in {\cal K}^{forbid}$. In practice, this set of rules is usually defined a priori by the operator to impose radiality constraints. Finally, constraints \eqref{eq_model4_upper_primal34} impose the binary nature of the switching variables. 

Following the decision-making process, the post-contingency operational problem is formulated in \eqref{eq_model4_lower_primal00}--\eqref{eq_model4_lower_primal31}:  
\begin{align}
    & H(\boldsymbol{z}^{sw},\boldsymbol{a}^L) = \underset{{\substack{
                            \Delta D^{p-^c}_b, \Delta D^{p+^c}_b, \\ 
                            \Delta D^{q-^c}_b, \Delta D^{q+^c}_b, \\
                            f^{p^c}_l, f^{q^c}_l, p^{{\color{black} {tr}^c}}_b,
                            q^{{\color{black}{tr}^c}}_b, v^{\dagger^c}_b                           
                           }}}{\text{Minimize}}
    ~ \sum_{b\in \mathcal{N}^{subs}}{C^{tr}_b} p^{{\color{black}{tr}}^c}_b \notag \\
    & + C^{\color{black}{ll}} \sum_{b\in \mathcal{N}}{\Big [ \Delta D^{p+^c}_b + \Delta D^{p-^c}_b + \Delta D^{q+^c}_b + \Delta D^{q-^c}_b\Big ]}  \label{eq_model4_lower_primal00} \\ \notag
    & ~ \text{subject to:} \notag \\ 
    & p^{{\color{black}{tr}^c}}_b + \sum_{l\in\mathcal{L}|to(l)=b}{f^{p^c}_l} - \sum_{l\in\mathcal{L}|fr(l)=b}{f^{p^c}_l} \notag \\
    & \hspace{15pt}  - D^p_b - \Delta D^{p+^c}_b + \Delta D^{p-^c}_b = 0 :(\eta_b^1); \forall b \in \mathcal{N}^{subs} \label{eq_model4_lower_primal01} \\
    & q^{{\color{black}{tr}^c}}_b + \sum_{l\in\mathcal{L}|to(l)=b}{f^{q^c}_l} - \sum_{l\in\mathcal{L}|fr(l)=b}{f^{q^c}_l} \notag \\
    & - \tan(\arccos(PF_b))D^p_b - \Delta  D^{q+^c}_b + \Delta  D^{q-^c}_b = 0: \notag \\
    & \hspace{155pt} (\eta_b^2); 
    \forall b \in \mathcal{N}^{subs} \label{eq_model4_lower_primal02} \\
    & \sum_{l\in\mathcal{L}|to(l)=b}{f^{p^c}_l} - \sum_{l\in\mathcal{L}|fr(l)=b}{f^{p^c}_l} - D^p_b - \Delta D^{p+^c}_b\notag \\
    & \hspace{64pt}  + \Delta D^{p-^c}_b = 0 :(\eta_b^3); \forall b \in \mathcal{N} \setminus \mathcal{N}^{subs} \label{eq_model4_lower_primal03} \\
    & \sum_{l\in\mathcal{L}|to(l)=b}{f^{q^c}} - \sum_{l\in\mathcal{L}|fr(l)=b}{f^{q^c}_l}  \notag \\
    & \hspace{10pt} - \tan(\arccos(PF_b))D^p_b - \Delta D^{q+^c}_b + \Delta D^{q-^c}_b = 0:  \notag \\
    & \hspace{136pt} (\eta_b^4); \forall b \in \mathcal{N} \setminus \mathcal{N}^{subs} \label{eq_model4_lower_primal04} \\
    & - v^{\dagger ^c}_{fr(l)} + v^{\dagger ^c}_{to(l)} + 2(R_l f^{p^c}_l + X_l f^{q^c}_l) \notag \\
    & \hspace{15pt} - (1 - a^L_l)M - (1 - z^{sw}_l)M \leq 0:(\eta_l^5); \forall l \in \mathcal{L}^{sw} \label{eq_model4_lower_primal05} \\
    & v^{\dagger ^c}_{fr(l)} - v^{\dagger ^c}_{to(l)} - 2(R_l f^{p^c}_l + X_l f^{q^c}_l) - (1 - a^L_l)M  \notag \\
    & \hspace{74pt} - (1 - z^{sw}_l)M \leq 0:(\eta_l^6); \forall l \in \mathcal{L}^{sw} \label{eq_model4_lower_primal06} \\
    & - v^{\dagger ^c}_{fr(l)} + v^{\dagger ^c}_{to(l)} + 2(R_l f^{p^c}_l + X_l f^{q^c}_l) \notag \\
    & \hspace{62pt} - (1 - a^L_l)M \leq 0:(\eta_l^7); \forall l \in \mathcal{L} \setminus \mathcal{L}^{sw} \label{eq_model4_lower_primal07} \\
    & v^{\dagger ^c}_{fr(l)} - v^{\dagger ^c}_{to(l)} - 2(R_l f^{p^c}_l + X_l f^{q^c}_l) \notag \\
    & \hspace{62pt} - (1 - a^L_l)M \leq 0: (\eta_l^8); \forall l \in \mathcal{L} \setminus \mathcal{L}^{sw} \label{eq_model4_lower_primal08} \\
    & \underline{V}_b^2 \leq v^{\dagger ^c}_{b} \leq \overline{V}_b^2:(\eta_b^9,\eta_b^{10}); \forall b \in \mathcal{N} \label{eq_model4_lower_primal09} \\
    %
    %
    & - z^{sw}_l \overline{F}_l \leq f^{p^c}_l \leq z^{sw}_l \overline{F}_l: (\eta_l^{11},\eta_l^{12}); \forall l \in {\mathcal L}^{sw} \label{eq_model4_lower_primal11} \\   
    %
    %
    & - z^{sw}_l \overline{F}_l \leq f^{q^c}_l \leq z^{sw}_l \overline{F}_l: (\eta_l^{13},\eta_l^{14}); \forall l \in {\mathcal L}^{sw} \label{eq_model4_lower_primal13} \\
    %
    %
    & - a^{L}_l \overline{F}_l \leq f^{p^c}_l \leq a^{L}_l \overline{F}_l:(\eta_l^{15},\eta_l^{16}); \forall l \in {\mathcal L} \label{eq_model4_lower_primal15} \\
    %
    %
    & - a^{L}_l \overline{F}_l \leq f^{q^c}_l \leq a^{L}_l \overline{F}_l:(\eta_l^{17},\eta_l^{18}); \forall l \in {\mathcal L} \label{eq_model4_lower_primal17} \\
    %
    %
    & f^{q^c}_l - {\color{black}\cot} \Biggl ( \Biggl ( \frac{1}{2} - e  \Biggr ) \frac{\pi}{4} \Biggr ) \Biggl ( f^{p^c}_l - {\color{black}\cos} \Biggl ( e \frac{\pi}{4} \Biggr ) \overline{F}_l \Biggr ) \notag \\
    & \hspace{14pt} - {\color{black}\sin} \Biggl ( e \frac{\pi}{4} \Biggr )\overline{F}_l \leq 0:(\eta_{l,e}^{19}); \forall l \in {\mathcal L}, e \in \{1,\dots,4\} \label{eq_model4_lower_primal19} \\
    & - f^{q^c}_l - {\color{black}\cot} \Biggl ( \Biggl ( \frac{1}{2} - e  \Biggr ) \frac{\pi}{4} \Biggr ) \Biggl ( f^{p^c}_l - {\color{black}\cos} \Biggl ( e \frac{\pi}{4} \Biggr ) \overline{F}_l \Biggr ) \notag \\
    & \hspace{14pt} - {\color{black}\sin} \Biggl ( e \frac{\pi}{4} \Biggr )\overline{F}_l \leq 0: (\eta_{l,e}^{20}); \forall l \in {\mathcal L}, e \in \{1,\dots,4\} \label{eq_model4_lower_primal20} \\
    & 0 \leq p^{{\color{black} {tr}^c }}_b \leq \overline{P}_b^{tr^c}: (\eta_b^{21},\eta_b^{22}); \forall b \in \mathcal{N}^{subs} \label{eq_model4_lower_primal21} \\
    %
    %
    & \underline{Q}^{tr^c}_b \leq q^{{\color{black} {tr}^c }}_b \leq \overline{Q}^{tr^c}_b:(\eta_b^{23},\eta_b^{24}); \forall b \in \mathcal{N}^{subs} \label{eq_model4_lower_primal23} \\
    %
    %
    & v^{\dagger^c}_b = V^{ref^2}:(\eta_b^{25}); \forall b \in \mathcal{N}^{subs} \label{eq_model4_lower_primal25} \\    
    & \Delta D_b^{p^{+^c}},\Delta D_b^{p^{-^c}},\Delta D_b^{q^{+^c}},\Delta D_b^{q^{-^c}} \geq 0: \notag \\
    & \hspace{110pt}(\eta_b^{26},\eta_b^{27},\eta_b^{28},\eta_b^{29}); \forall b \in \mathcal{N} \label{eq_model4_lower_primal26} \\
    %
    %
    %
    %
    & \Delta D_b^{p^{-^c}} \leq D^p_b: (\eta_b^{30}); \forall b \in \mathcal{N} \label{eq_model4_lower_primal30} \\
    & \Delta D_b^{q^{-^c}} \leq \tan(\arccos(PF_b)) D^p_b: (\eta_b^{31}); \forall b \in \mathcal{N} \label{eq_model4_lower_primal31}
\end{align}


\noindent where the symbols within parenthesis are the dual variables associated with the constraints. Problem \eqref{eq_model4_lower_primal00}--\eqref{eq_model4_lower_primal31} is a linear programming problem with (continuous) decision variables $p^{{\color{black}{tr}}^c}_{b}$, $q^{{\color{black}{tr}}^c}_{b}$, ${\color{black}v^{\dagger^c}_{b}}$, $f^{p^c}_{l}$, $f^{q^c}_{l}$, ${\Delta}D^{p+^c}_{b}$, ${\Delta}D^{p-^c}_{b}$, ${\Delta}D^{q+^c}_{b}$, and ${\Delta}D^{q-^c}_{b}$ with essentially the same role as in \eqref{eq_model4_upper_primal00}--\eqref{eq_model4_upper_primal34}. Analogously to \eqref{eq_model4_upper_primal01}--\eqref{eq_model4_upper_primal04}, constraints \eqref{eq_model4_lower_primal01}--\eqref{eq_model4_lower_primal04} model active and reactive power balances. Constraints \eqref{eq_model4_lower_primal05}--\eqref{eq_model4_lower_primal08} express voltage differences in line segments under a given contingency state associated with vector $\boldsymbol{a}^{L}$ and a first-stage switching decision $z^{sw}_l$. Constraints \eqref{eq_model4_lower_primal09} impose voltage limits. Constraints \eqref{eq_model4_lower_primal11}--\eqref{eq_model4_lower_primal20} enforce limits to active and reactive flows. 
Constraints \eqref{eq_model4_lower_primal21}--\eqref{eq_model4_lower_primal31} limit power injections and impose voltage reference at the substations as well as enforce non-negativity to power surplus and power loss variables.


\subsection{Decision-(Line-Flows)-Dependent Ambiguity Set Modeling} \label{sec::AmbiguitySet}

Following the discussion of the previous section, the proposed methodology for distribution system operations seeks for least-cost pre- and post-contingency states operative decisions, the latter with respect to line segment availability. We argue, furthermore, that such line availability is mainly influenced by exogenous weather conditions, in particular during adverse climate circumstances, as well as endogenously impacted by the determined operative point and power flow in the network \cite{Muhs2020}. To jointly tackle these two critical uncertain factors in a unified framework, in this section, a pre-contingency line-flow-dependent ambiguity set of credible branch availability probabilities is constructed. More specifically, the uncertainty related to the underlying stochastic process associated with line failures is modeled via a tailored ambiguity set ${\cal P} \in {\cal M}_{+}$ composed of a collection of probability distributions that characterize the limited knowledge of failure probabilities and the endogenous/exogenous uncertain impact factors. Formally, the proposed ambiguity set is expressed as:
\begin{align}
    & {\cal P}(\boldsymbol{f}^{p}, \boldsymbol{\beta}) = \Bigl\{ {\cal Q} \in {\cal M}_{+} ({\cal A}) ~ \Big| ~ \mathbb{E}_{\cal Q} \bigl[\boldsymbol{S} \hat{\boldsymbol{a}}^{L}\big] \leq \overline{\boldsymbol{\mu}}(\boldsymbol{f}^{p}, {\boldsymbol{\beta}}) \Bigr\}. \label{AmbiguitySet_1}
\end{align}

In \eqref{AmbiguitySet_1}, function $\overline{\boldsymbol{\mu}} (\cdot, \cdot)$ is a vector of means that defines the dependency of external factors and decisions variables. The term $\boldsymbol{S}$ is defined as an auxiliary matrix of coefficients, and $\hat{\boldsymbol{a}}^{L} = \mathbb{1} - \boldsymbol{a}^{L}$ indicates a random vector of \textit{line unavailability} with set ${\cal A}$ characterizing its support. In this work, the support of the random vector $\boldsymbol{a}^{L}$, is defined as
%
\begin{align}
	& {\cal A}=\biggl \{\boldsymbol{a}^{L} \in \{0,1\}^{|{\cal L}|} ~ \bigg| ~ \sum_{l \in {\cal L}} a^L_l \geq |{\cal L}| - K \biggr \}, \label{AmbiguitySet_2}
\end{align}

\noindent with $K$ indicating the number of simultaneous unavailable system components (lines segments, in the context of this work) \cite{Moreira2018_EnergyReserveSchedu_AmbiRen, Moreira2019_AmbiguityTEP}. Following the ambiguity set definition \eqref{AmbiguitySet_1}, fundamentally, a critical modeling element is the appropriate definition of the vector of means $\bigl(\overline{\boldsymbol{\mu}}(\boldsymbol{f}^{p}, {\boldsymbol{\beta}})\bigr)$. In this work, we follow the main findings in \cite{Muhs2020} and consider the following functional representation:
\begin{align}
    \overline{\boldsymbol{\mu}}(\boldsymbol{f}^{p}, {\boldsymbol{\beta}}) = \boldsymbol{\gamma} + {\color{black}\mathrm{diag}}(\boldsymbol{\beta}) |\boldsymbol{f}^p|, \label{AmbiguitySet_3}
\end{align}

\noindent where ${\color{black}\mathrm{diag}}(\boldsymbol{\beta})$ returns a diagonal matrix with elements of $\boldsymbol{\beta}$. Structurally, vector $\boldsymbol{\gamma}$ represents an estimated upper bound for the nominal probability of failure associated with each line segment $l \in {\cal L}$, extracted from the set of available information (e.g., failures per year), whereas vector $\boldsymbol{\beta}$ (exogenous-impact) characterizes the sensitivity in the probability of failure to the scheduled active power flow  (endogenous-impact) in each line. Within the context of this paper, on the one hand, vector $\boldsymbol{\beta}$ provides instrumental information on how the probability of line failure increases as a function of the power flows. On the other hand, in particular, during adverse climate conditions (e.g., dry weather and high wind speed), the line failure can be caused by fire, started by the line itself if it is sufficiently close to vegetation. This condition can be adjusted by the system operator using the contextual (exogenous) vector $\boldsymbol{\beta}$. Therefore, structurally, by setting $\boldsymbol{S} = \big[\mathbb{I} ~ | ~ -\mathbb{I}\big]^T_{2|{\cal L}| \times |{\cal L}|}$ and $\overline{\mu}_{l} = (\gamma_l + \beta_l |f^p_{l}|), ~ \forall ~ l \in {\cal L}$ and $\overline{\mu}_{(l+|{\cal L}|)} = 0, ~ \forall ~ l \in {\cal L}$ in \eqref{AmbiguitySet_3}, we have the resulting ambiguity set:
\begin{align}
    & \hspace{-0.00cm} {\cal P}(\boldsymbol{f}^{p}, \boldsymbol{\beta}) = \Bigl\{ {\cal Q} \in {\cal M}_+({\cal A}) ~ \Big| ~ 0 \leq \mathbb{E}_{\cal Q} [\hat{a}^L_l] \leq \gamma_{l} + \beta_{l} |{f}^{p}_{l}|; \notag \\
    & \hspace{6.20cm} \forall ~ l \in {\cal L} \Bigr\}. \label{AmbiguitySet_4}
\end{align}

\noindent Since $\hat{\boldsymbol{a}}^{L} = \mathbb{1} - \boldsymbol{a}^{L}$ is a Bernoulli-type random vector, the structural specification of \eqref{AmbiguitySet_4} implies that a failure probability in each line $l \in {\cal L}$ is constrained by the factor $\gamma_l + \beta_l |{f}^p_l|$, thus dependent on the (endogenous) scheduled active power flow ${f}^{p}_{l}$ and the contextual (exogenous) information $\beta_{l}$. It is worth highlighting that the proposed distribution system operation model \eqref{eq_model4_upper_primal00}--\eqref{eq_model4_upper_primal34} with \eqref{AmbiguitySet_4} has a decision process that follows a two-stage, distributionally robust optimization with decision-dependent ambiguity set rationale. This decision process is formulated as a three-level system of optimization problems, not suitable for direct implementation on commercial solvers nor standard mathematical programming algorithms. Therefore, in the next section, we leverage the problem structure to devise a decomposition-based solution approach to efficiently handle the proposed model.

\section{Solution methodology}\label{SolutionMethodology}

The two-stage formulation \eqref{eq_model4_upper_primal00}--\eqref{eq_model4_upper_primal34} proposed in Section \ref{MathematicalFormulation} is intended to model the operation of a distribution system while performing switching actions to minimize the worst-case expected cost in post-contingency operations. In this section, we develop an iterative procedure based on outer approximation to solve this problem. We being by replacing the last term in \eqref{eq_model4_upper_primal00} with $\alpha$ and writing \eqref{Model_3_1}. The variable $\alpha$ is defined through \eqref{Model_3_3}--\eqref{Model_3_5} which essentially represents the last term in \eqref{eq_model4_upper_primal00}. Thus, we equivalently rewrite model \eqref{eq_model4_upper_primal00}--\eqref{eq_model4_upper_primal34}  as \eqref{Model_3_1}--\eqref{Model_3_5}.
%
%
\begin{align}
	& \underset{{\substack{ 
            \alpha, {\Delta}D^{p-}_{b}, {\Delta}D^{p+}_{b}, {\Delta}D^{q-}_{b}, {\Delta}D^{q+}_{b}, \\
            f^p_{l}, f^q_{l}, p^{{\color{black} tr }}_{b}, q^{{\color{black} tr }}_{b}, v^{\dagger}_{b}, y^{sw}_{l}, z^{sw}_{l}}}}{\text{Minimize}}  \hspace{0.1cm} \sum_{b \in {\color{black}\mathcal{N}^{subs}}} C^{tr}_{b} p^{{\color{black} tr }}_{b} \notag\\
	& + \sum_{b \in {\color{black}\mathcal{N}}} C^{\color{black}{ll}}\Bigl({\Delta}D^{p+}_{b} + {\Delta}D^{p-}_{b} + {\Delta}D^{q+}_{b} + {\Delta}D^{q-}_{b}\Bigr)  \notag\\
	& + \sum_{l \in {\cal L}^{sw}} C^{sw}_l  y^{sw}_{l} + \alpha \label{Model_3_1}\\
	& \text{subject to:}\notag\\
	& \text{Constraints \eqref{eq_model4_upper_primal01}--\eqref{eq_model4_upper_primal34}} \label{Model_3_2}\\
	& \alpha = \biggl \{ \underset{{\substack{ {\cal Q} \in {\cal M}_{+} }}}{\text{Maximize}} \hspace{0.1cm} \sum_{\boldsymbol{a}^L \in {\cal A}} H\bigl(\boldsymbol{z}^{sw},\boldsymbol{a}^L\bigr) {\cal Q}(\boldsymbol{a}^L) \label{Model_3_3}\\
	& \hspace{28pt} \text{subject to:}\notag\\
	& \hspace{28pt} \sum_{\boldsymbol{a}^L \in {\cal A}} \bigl(\boldsymbol{S} \hat{\boldsymbol{a}}^L\bigr) {\cal Q}(\boldsymbol{a}^L) \leq \overline{\boldsymbol{\mu}}(\boldsymbol{f}^{p}, {\boldsymbol{\beta}}) : (\boldsymbol{\psi}) \label{Model_3_4}\\
	& \hspace{28pt} \sum_{\boldsymbol{a}^L \in {\cal A}} {\cal Q}(\boldsymbol{a}^L) = 1 : ({\varphi}) \biggr\} 
	\label{Model_3_5}.
\end{align}
Resorting to duality theory, we can substitute $\alpha$ in \eqref{Model_3_1} by the dual objective function of the inner model \eqref{Model_3_3}--\eqref{Model_3_5} and replace it with the dual feasibility constraints. More precisely,
\begin{align}
	& \underset{{\substack{
                          {\Delta}D^{p-}_{b}, {\Delta}D^{p+}_{b}, {\Delta}D^{q-}_{b}, {\Delta}D^{q+}_{b}, \varphi, \\
                           \boldsymbol{\psi} \geq \boldsymbol{0}, f^p_{l}, f^q_{l}, p^{{\color{black} tr }}_{b}, q^{{\color{black} tr}}_{b}, {\color{black}v^{\dagger}_{b}}, y^{sw}_{l}, z^{sw}_{l}}}}{\text{Minimize}} \hspace{0.01cm} \sum_{b \in {\color{black}\mathcal{N}^{subs}}} C^{tr}_{b} p^{{\color{black} tr }}_{b} \notag\\
	& \hspace{0.2cm} + \sum_{b \in {\color{black}\mathcal{N}}} C^{\color{black}{ll}}\Bigl({\Delta}D^{p+}_{b}+ {\Delta}D^{p-}_{b} + {\Delta}D^{q+}_{b} + {\Delta}D^{q-}_{b}\Bigr)  \notag\\
	& \hspace{0.2cm} + \sum_{l \in {\cal L}^{sw}} C^{sw}_l  y^{sw}_{l} + \boldsymbol{\psi}^{\top} \overline{\boldsymbol{\mu}}(\boldsymbol{f}^{p}, {\boldsymbol{\beta}}) + \varphi \label{Model_4_1}\\
	& \text{subject to:} \notag\\
	& \text{Constraints \eqref{eq_model4_upper_primal01}--\eqref{eq_model4_upper_primal34}} \label{Model_4_2}\\
	& \boldsymbol{\psi}^{\top} \boldsymbol{S} \hat{\boldsymbol{a}}^L + \varphi \geq H\bigl(\boldsymbol{z}^{sw},\boldsymbol{a}^L\bigr); \forall ~ \boldsymbol{a}^L \in {\cal A} \label{Model_4_3}.
\end{align}

\noindent To withstand the intractability caused by the combinatorial nature of the support set ${\cal A}$ defined in \eqref{AmbiguitySet_2}, we replace constraints in \eqref{Model_4_3} by:
\begin{align}
    & \varphi \geq \max_{\boldsymbol{a}^L \in {\cal A}} \Bigl \{ H\bigl(\boldsymbol{z}^{sw},\boldsymbol{a}^L\bigr) - \boldsymbol{\psi}^{\top} \boldsymbol{S} \hat{\boldsymbol{a}}^L \Bigr \}
    \label{Model_4_4}.
\end{align}

Based on \eqref{Model_4_1}, \eqref{Model_4_2}, \eqref{Model_4_4}, we propose in the next subsections an iterative procedure to address formulation \eqref{eq_model4_upper_primal00}--\eqref{eq_model4_upper_primal34}. 

\vspace{-0.3cm}

\subsection{Subproblem} \label{subproblem}


The role of the subproblem is to provide an approximation to the right-hand side of \eqref{Model_4_4}. Note that $H(\boldsymbol{z}^{sw},\boldsymbol{a}^L)$ is a minimization problem. Thus, to build the subproblem, we take the following steps: (i) write the dual problem of $H(\boldsymbol{z}^{sw},\boldsymbol{a}^L)$, (ii) subtract the dual objective function by $\boldsymbol{\psi}^{\top} \boldsymbol{S} \hat{\boldsymbol{a}}^L$, and (iii) handle the bilinear products between dual and binary variables $\boldsymbol{a}^{L}$ in the dual objective function. It is worth mentioning that the recourse function associated with the resulting subproblem is convex with respect to the first-stage decision as it is a maximum of affine functions, therefore rendering the description of the right-hand side of \eqref{Model_4_4} suitable to cutting planes approximation.

\subsection{Master problem}
\vspace{-0.1cm}
The master problem developed in this Section is a relaxation of the original model \eqref{eq_model4_upper_primal00}--\eqref{eq_model4_upper_primal34}. Such relaxation is improved by the iterative inclusion of cutting planes. The master problem is formulated as follows.
\begin{align}
	& \underset{{\substack{ 
            {\Delta}D^{p-}_{b}, {\Delta}D^{p+}_{b}, {\Delta}D^{q-}_{b}, {\Delta}D^{q+}_{b}, \\
            {\color{black}\delta_{le}}, \xi_{l}, \rho_{le}, \varphi, \chi_{l}, \boldsymbol{\psi} \geq \boldsymbol{0}, f^p_{l}, \\
            f^{p,-}_{l}, f^{p,+}_{l}, f^q_{l}, p^{{\color{black} tr }}_{b}, q^{{\color{black} tr }}_{b}, {\color{black}v^{\dagger}_{b}}, y^{sw}_{l}, z^{sw}_{l}}}}{\text{Minimize}} \hspace{0.1cm} \sum_{b \in {\color{black}\mathcal{N}^{subs}}} C^{tr}_{b} p^{{\color{black} tr }}_{b} \notag\\
	& \hspace{0.50cm} + \sum_{b \in {\color{black}\mathcal{N}}} C^{\color{black}{ll}}\Bigl({\Delta}D^{p+}_{b}+ {\Delta}D^{p-}_{b} + {\Delta}D^{q+}_{b} + {\Delta}D^{q-}_{b}\Bigr) \notag \\
	& \hspace{0.50cm} + \sum_{l \in {\cal L}^{sw}} C^{sw}_l  y^{sw}_{l} + \sum_{l \in {\cal L}}  ({\color{black}\gamma_l} \psi_{l}  + \beta_l \chi_{l}) + \varphi \label{Model_9_1}\\
	& \text{subject to:} \notag \\
	& \text{Constraints \eqref{eq_model4_upper_primal01}--\eqref{eq_model4_upper_primal34}} \label{Model_9_2}\\
	& f^p_{l} = f^{p,+}_{l} - f^{p,-}_{l}; \forall l \in {\cal L} \label{Model_9_3} \\
	& 0 \leq f^{p,+}_{l} \leq {\color{black}\overline{F}}_l \xi_{l}; \forall l \in {\cal L} \label{Model_9_4} \\
	& 0 \leq f^{p,-}_{l} \leq {\color{black}\overline{F}}_l(1-\xi_{l}); \> \forall l \in {\cal L} \label{Model_9_5} \\
    & \xi_{l} \in \{0,1\}; \forall l \in {\cal L}\label{Model_9_6} \\
	& f^{p,+}_{l} + f^{p,-}_{l} = s \sum_{{\color{black}e=1}}^{E_l} 2^{{\color{black}e-1}} {\color{black}\delta_{le}}; \> \forall ~ l \in {\cal L} \label{Model_9_7}\\
	& \delta_{le} \in \{0,1\}; \forall l \in {\cal L}, e \in 1,\dots, E_l \label{Model_9_8} \\
	& -M(1 - \delta_{le}) \leq \psi_{l} - \rho_{le} \leq M(1 - \delta_{le}); \forall l \in {\cal L}, \notag \\
    & \hspace{132pt} e = 1,\dots, E_l \label{Model_9_9} \\
    & - \delta_{le} M \leq \rho_{le} \leq  \delta_{le} M; \> \forall l \in {\cal L}, e = 1,\dots, E_l \label{Model_9_10} \\
    & \chi_{l} = s \sum_{e = 1}^{E_l} 2^{e-1} \rho_{le}; \forall l \in {\cal L} \label{Model_9_11} \\
    & \varphi \geq \sum_{b \in \mathcal{N}^{subs}} \Biggl [ - D^p_b \eta_b^{1^{(j)}} -\tan{(\arccos{(PF_b)})} D^p_b \eta_b^{2^{(j)}} \notag \\
    & \hspace{2pt} +\underline{V}^2_{b} \eta_b^{9^{(j)}} - \overline{V}^2_{b} \eta_b^{10^{(j)}} - \overline{P}^{tr^c}_b \eta_b^{22^{(j)}} +  \underline{Q}^{tr^c}_b \eta_b^{23^{(j)}} \notag \\
    & \hspace{2pt} - \overline{Q}^{tr^c}_b \eta_b^{24^{(j)}} - V^{ref^2} \eta_b^{25^{(j)}} - D^p_b  \eta_b^{30^{(j)}} \notag \\
    & \hspace{2pt} - \tan(\arccos(PF_b)) D^p_b  \eta_b^{31^{(j)}} \Biggr ] \notag \\
    & \hspace{2pt} + \sum_{b \in \mathcal{N} \setminus \mathcal{N}^{subs}} \Biggl [ - D^p_b \eta_b^{3^{(j)}} -\tan{(\arccos{(PF_b)})} D^p_b \eta_b^{4^{(j)}} \notag \\
    & \hspace{2pt} + \underline{V}^2_{b} \eta_b^{9^{(j)}} - \overline{V}^2_{b} \eta_b^{10^{(j)}} - D^p_b  \eta_b^{30^{(j)}} \notag \\
    & \hspace{2pt} - \tan(\arccos(PF_b)) D^p_b  \eta_b^{31^{(j)}} \Biggr ] \notag \\
    & \hspace{2pt} + \sum_{l \in \mathcal{L} \setminus \mathcal{L}^{sw}} \Biggl [ - (1 - a^{L^{(j)}}_l) M \eta_l^{7^{(j)}} - (1 - a^{L^{(j)}}_l) M \eta_l^{8^{(j)}} \notag \\
    & \hspace{2pt} - a^{L^{(j)}}_l \overline{F}_l (\eta_l^{15^{(j)}} + \eta_l^{16^{(j)}} + \eta_l^{17^{(j)}} + \eta_l^{18^{(j)}}) \notag \\
    & \hspace{2pt} + \sum_{e \in \{1,2,3,4\}} \Bigg ( \overline{F}_l \Big( {\color{black}\cot} \Big(((1/2) - e) (\pi/4)\Big) {\color{black}\cos} \Big( e(\pi/4) \Big) \notag \\
    & \hspace{2pt} - {\color{black}\sin} \Bigr ( e(\pi/4) \Bigr ) \Bigr ) (\eta_{l,e}^{19^{(j)}} + \eta_{l,e}^{20^{(j)}}) \Biggr ) \Biggr ] \notag \\
    & \hspace{2pt} + \sum_{l \in \mathcal{L}^{sw}} \Bigg [ - ( (1 - a^{L^{(j)}}_l) M + (1 - z^{sw}_l) M ) \eta_l^{5^{(j)}} \notag \\
    & \hspace{2pt} - ( (1 - a^{L^{(j)}}_l) M + (1 - z^{sw}_l) M ) \eta_l^{6^{(j)}} - z^{sw}_l \overline{F}_l (\eta_l^{11^{(j)}} \notag \\
    & \hspace{2pt} + \eta_l^{12^{(j)}} + \eta_l^{13^{(j)}} + \eta_l^{14^{(j)}}) - a^{L^{(j)}}_l \overline{F}_l (\eta_l^{15^{(j)}} + \eta_l^{16^{(j)}} \notag \\
    & \hspace{2pt} + \eta_l^{17^{(j)}} + \eta_l^{18^{(j)}}) \notag \\
    & \hspace{2pt} + \sum_{e \in \{1,2,3,4\}} \Bigg ( \overline{F}_l \Big( {\color{black}\cot} \Big(((1/2) - e) (\pi/4)\Big) {\color{black}\cos} \Big( e(\pi/4) \Big) \notag \\
    & \hspace{2pt} - {\color{black}\sin} \Big( e(\pi/4) \Big)\Big) (\eta_{l,e}^{19^{(j)}} + \eta_{l,e}^{20^{(j)}}) \Bigg ) \Bigg ] - \notag \\
    & \hspace{2pt} \sum_{l \in \mathcal{L}} \Bigg [ (\psi_l - \psi_{|\mathcal{L}|+l})(1 - a^{L^{(j)}}_l) \Bigg ] \forall ~ j \in {\cal J} \label{Model_9_12},
\end{align}

\noindent where the product $\boldsymbol{\psi}^{\top} \overline{\boldsymbol{\mu}}(\boldsymbol{f}^{p}, {\boldsymbol{\beta}})$ in the objective function is replaced by $\sum_{l \in {\cal L}}  ({\color{black}\gamma_l} \psi_{l}  + \beta_l \chi_{l})$ and $\chi_{l}$ represents the bilinear term ${\psi}_l |f^p_l|$ as modeled in \eqref{Model_9_3}--\eqref{Model_9_11}. Furthermore, expression \eqref{Model_9_12} represents cutting planes that are iteratively included to approximate the right-hand side of expression \eqref{Model_4_4}. 

\subsection{Solution Algorithm} \label{Sec::SolAlg}

In this section, we describe the outer approximation algorithm proposed in this work, following the Master and Subproblem descriptions. Structurally, it is an iterative process that is carried out until the approximation provided by the inclusion of the cutting planes \eqref{Model_9_12} is sufficient to make the solution of the relaxed Master problem close enough to optimality. This proposed outer approximation algorithm is summarized as follows.

\begin{enumerate}
	\item {Initialization:} set counter $m \leftarrow 0$ and set $\mathcal{J} \leftarrow \emptyset$. 
	\item Solve the optimization model \eqref{Model_9_1}--\eqref{Model_9_12}, store $\boldsymbol{z}^{sw (m)}$, $\boldsymbol{\psi}^{(m)}$ and $\varphi^{(m)}$, and set ${LB}^{(m)}$ equal to the value of the objective function \eqref{Model_9_1}.
	\item Identify the worst case contingency for $\boldsymbol{z}^{sw (m)}$ and $\boldsymbol{\psi}^{(m)}$ by running the linearized subproblem described in Subsection \ref{subproblem}. Store values of its decision variables and calculate ${UB}^{(m)}$ by subtracting $\varphi^{(m)}$ from ${LB}^{(m)}$ and adding the value of the objective function of the subproblem.
	\item If $\bigl(UB^{(m)} - LB^{(m)}\bigr)/UB^{(m)} \leq \epsilon$, then STOP; else, CONTINUE.
	\item Include in \eqref{Model_9_1}--\eqref{Model_9_12} a new cutting plane of the format \eqref{Model_9_12} with decision variables stored in Step 3, set $m \leftarrow m+1$, ${\cal J} \leftarrow {\cal J} \cup \{m\}$, and go to Step 2.
\end{enumerate}

It is interesting to note that the cuts generated when considering $\boldsymbol{\beta} = \boldsymbol{0}$ (i.e., neglecting decision-dependent uncertainty) are still valid for the decision-dependent case $(\boldsymbol{\beta} \geq \boldsymbol{0})$. This happens because the vectors of decision variables $\boldsymbol{z}^{sw}$ and $\boldsymbol{\psi}$ have the same feasible region regardless of the value of $\boldsymbol{\beta}$ and the cuts obtained by solving the maximization problem on the right-hand side of \eqref{Model_4_4} would be valid even if $\boldsymbol{z}^{sw}$ and $\boldsymbol{\psi}$ are not optimally decided by the Master problem. In the numerical experiments conducted in this work, we will leverage this property to accelerate the solution of the cases with decision-dependent uncertainty $(\boldsymbol{\beta} \geq \boldsymbol{0})$ by reusing the cutting planes obtained for the case where decision-dependent uncertainty is not considered ($\boldsymbol{\beta} = \boldsymbol{0}$). This reuse can be particularly advantageous since (i) it is usually much faster to solve the problem with $\boldsymbol{\beta} = \boldsymbol{0}$ and (ii) warming up the problem for $\boldsymbol{\beta} \geq \boldsymbol{0}$ with previously identified valid cutting planes can significantly improve computational efficiency as will be seen in the numerical experiments.

\begin{figure*}[!tb]
    \centering
     \includegraphics[width=.73\textwidth,height=0.5\textheight,keepaspectratio]{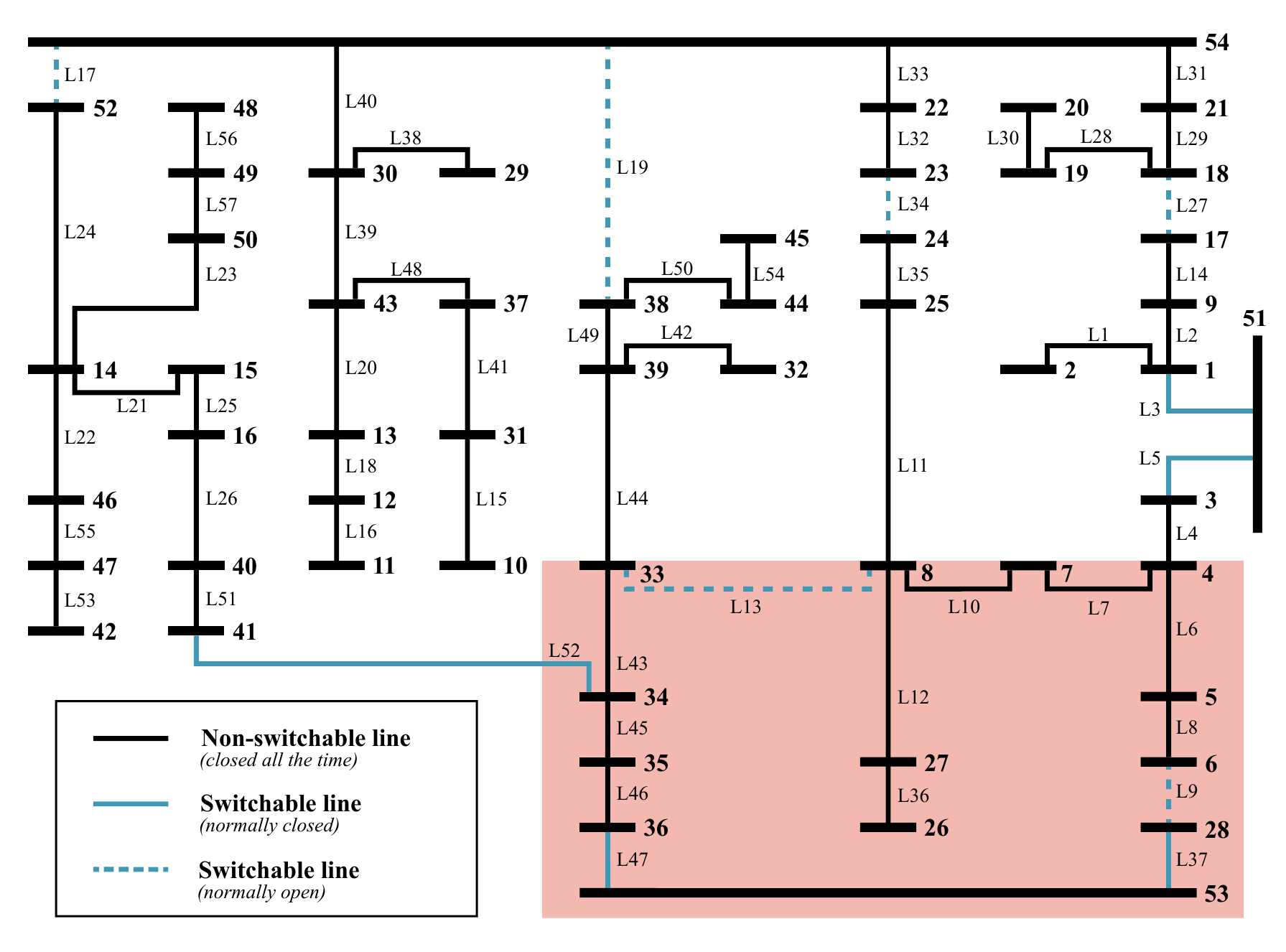}
    \caption{54-Bus distribution system.}
    \label{Fig:54BusFigure}
\end{figure*}

\begin{figure}[!tb]
    \centering
     \includegraphics[width=.45\textwidth,height=0.5\textheight,keepaspectratio]{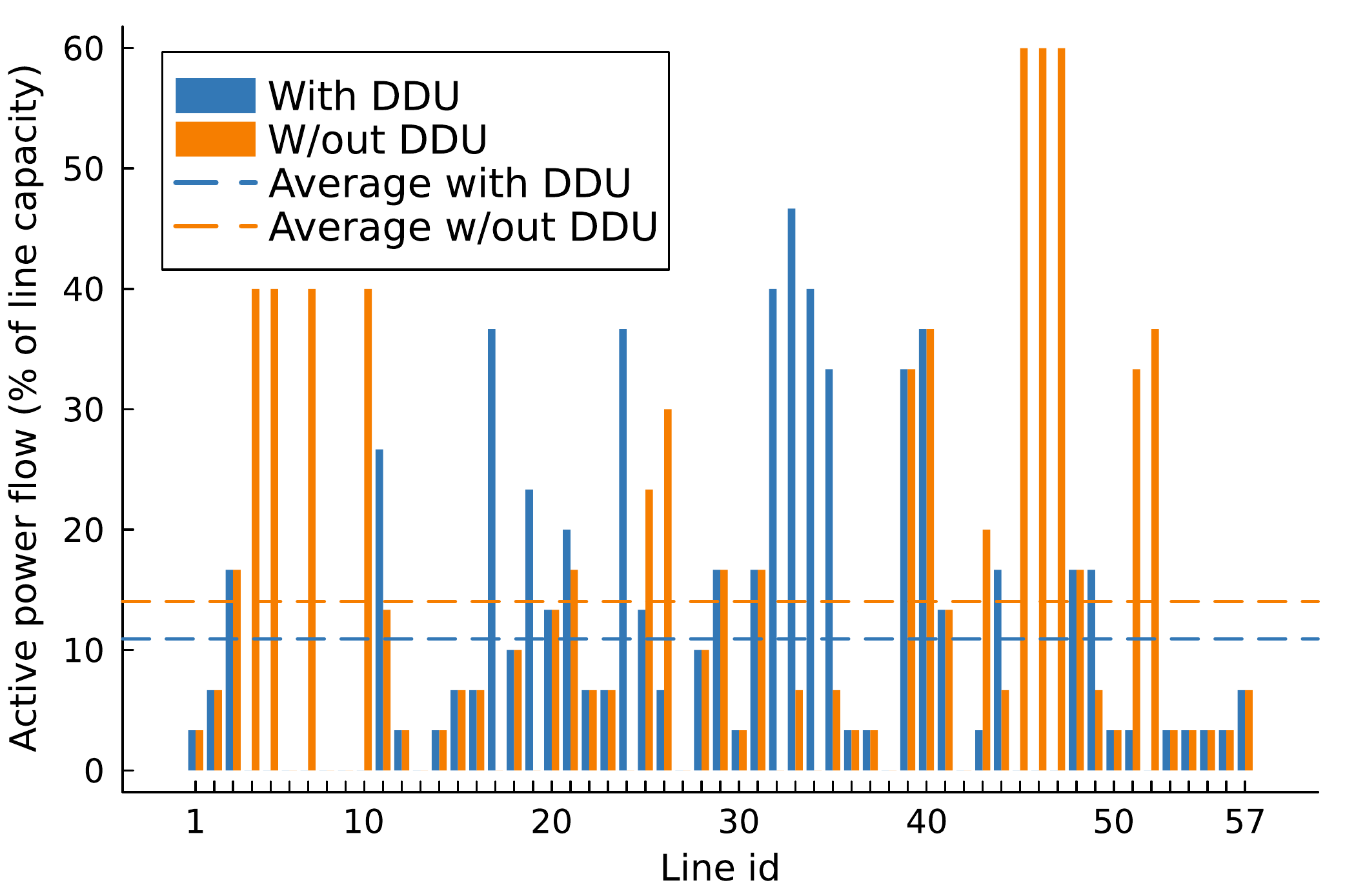}
    \caption{Power flows for the solutions \textit{with} and \textit{without DDU}.}
    \label{Fig:flowsComparison}
\end{figure}

\begin{figure}[!tb]
    \centering
     \includegraphics[width=.45\textwidth,height=0.5\textheight,keepaspectratio]{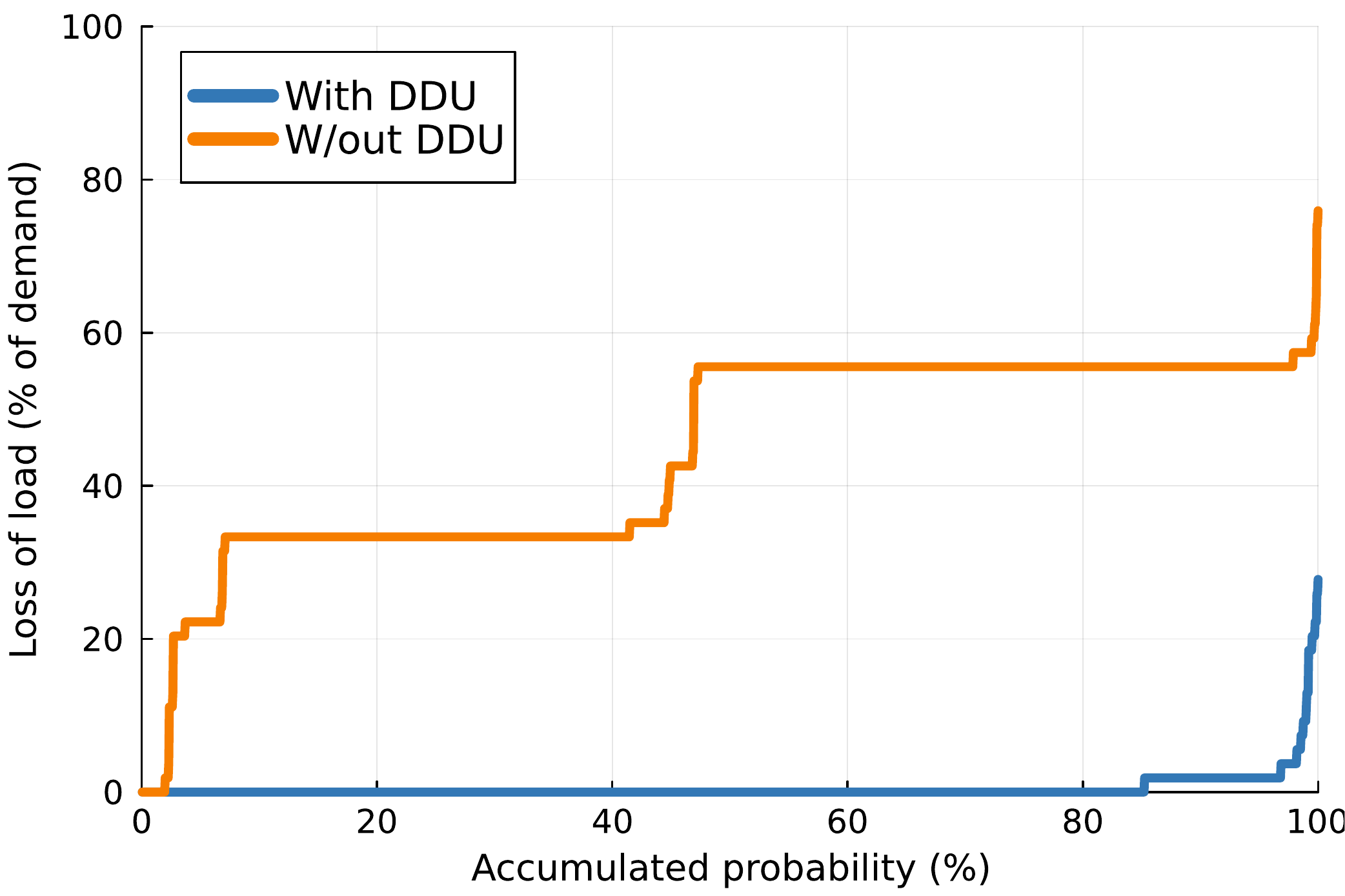}
    \caption{Out-of-sample inverse cumulative distribution of the system loss of load for the solutions \textit{with} and \textit{without DDU}.}
    \label{Fig:inverseCumulativeDist}
\end{figure}

\section{Case studies}\label{Sec:caseStudy}

The proposed methodology is illustrated in this section with two case studies. The first case study is based on a 54-bus distribution system, whereas the second one comprises a 138-bus distribution system. In both case studies, we consider that part of the grid is vulnerable to the ignition of a wildfire, which can be influenced by the levels of power flows passing through the line segments within the region. The solution algorithm described in Section \ref{Sec::SolAlg} has been implemented in Julia 1.6 and solved on a server with one Intel\textregistered ~Core\textregistered ~i7-10700K processor @ 3.80GHz and 64 GB of RAM, using Gurobi 9.0.3. under JuMP.

\subsection{54-bus system} \label{Sec:54bus}


In this case, we consider a 54-bus distribution system (depicted in Fig. \ref{Fig:54BusFigure}) based on the data provided in \cite{Munoz2016}. In this system, there are 3 substations (buses 51, 53, and 54 in Fig. \ref{Fig:54BusFigure}) and 57 lines. The total demand of the system is 5400 kW and the energy price is 0.01 \$/kWh. In addition, we consider that each switching action costs \$100, which can be performed in 11 out of the 57 lines. In Fig. \ref{Fig:54BusFigure}, the switchable lines are represented by blue lines. Furthermore, the blue dashed lines are initially open lines whereas the blue solid lines are initially closed. To enforce radiality constraints, lines L9 and L37 cannot be switched on simultaneously. The same rule applies to the pairs of lines L17 and L52, L13 and L47, L5 and L34, L3 and L27, L13 and L19, L19 and L47, and L13 and L47. 
These rules constitute the forbidden switching patterns in this case study. For replicability purposes, input data can be downloaded from \cite{CaseStudyRepository_DDU_wildfire2023}. In this case study, we consider an event of adverse climate conditions approaching that includes extreme dry weather and consistent wind speed. In addition, part of the grid, more specifically the southeast, is located close to vegetation, which renders this area particularly more likely to initiate a wildfire. The southeast area of the grid includes lines L6, L7, L8, L9, L10, L12, L13, L36, L37, L43, L45, L46, L47, and L52. We consider that every line segment has a nominal {\color{black}rate of failure equal to 0.4 failures per year}. Using the exponential probability distribution, this rate of failure translates into a failure probability of {\color{black}0.11\% for each line in the next 24 hours}. In addition, due to the adverse climate conditions, each of the aforementioned lines that belong to the southeast area has an increase of 3\% in its probability of failure for each 0.01 pu (100kW) of scheduled active power flow
. The remaining lines have an increase of $10^{-4}$\% in their probabilities of failure per 0.01 pu of scheduled active power flow.

Within this context, we consider three possible modeling and algorithmic structures to determine the status of switchable lines. In the first one, hereinafter referred to as {\it without DDU} (Decision-Dependent Uncertainty), the operator ignores the decision-dependent influence of line flows and probabilities of failures in the modeling and, therefore, only considers the nominal probabilities previously described. To do so, equation \eqref{AmbiguitySet_3} is modified to $\overline{\boldsymbol{\mu}} = \boldsymbol{\gamma}$. In the second one, hereinafter referred to as {\it with DDU}, the operator explicitly considers the aforementioned increase in failure probability corresponding to line usage according to \eqref{AmbiguitySet_3}. In the third one, hereinafter referred to as {\it with DDU and warm up}, decision-dependent is considered exactly as in the {\it with DDU} case but the cutting planes of the {\it without DDU} case are included in the master problem since the beginning of the execution of the solution algorithm. This reuse of cutting planes can help the {\it with DDU and warm up} approach to achieve the same solution of the {\it with DDU} method in less time. The respective switching statuses are depicted in Table \ref{tab:switchingSecisions}, where 1 means closed line and 0 means open line. As expected, when DDU is ignored, there is no incentive to change the status of any line since the nominal probabilities of failure are relatively low. Nonetheless, when DDU is considered, six lines have their statuses changed. In this context, the solution {\it without DDU} costs \$54, which is equivalent to the cost of feeding the loads without any switching, and the solution {\it with DDU} costs \$654, which includes feeding loads and performing 6 switching actions. The {\it with DDU and warm up} solution results in exactly the same costs and switching decisions as the {\it with DDU} solution. In Fig. \ref{Fig:flowsComparison}, it can be noted that the average flow per line, as well as the maximum flow among all branches, are significantly reduced when DDU is taken into account to decrease failure probabilities. {\color{black}The solutions {\it without DDU}, {\it with DDU}, and {\it with DDU and warm up}, were obtained in 10.59s, 49.22s, and 22.50s, respectively.}

\begin{table}[htbp]
    {\color{black}
    \centering
    \caption{Switching decisions with and without DDU}
    \renewcommand{\arraystretch}{1.30}
    \resizebox{\columnwidth}{!}{%
    \begin{tabular}{r | c c c c c c c c c c c |}
                         & \multicolumn{11}{c|}{\textbf{Switchable Lines}}                                                                                                       \\ \cline{2-12} 
                         & \textbf{3} & \textbf{5} & \textbf{9} & \textbf{13} & \textbf{17} & \textbf{19} & \textbf{27} & \textbf{34} & \textbf{37} & \textbf{47} & \textbf{52} \\ \cline{1-12} 
    \textbf{W/out DDU} & 1          & 1          & 0          & 0           & 0           & 0           & 0           & 0           & 1           & 1           & 1           \\ 
    \textbf{With DDU}    & 1          & 0          & 0          & 0           & 1           & 1           & 0           & 1           & 1           & 0           & 0          \\
    \textbf{{\color{black}With DDU and warm up}}    & 1          & 0          & 0          & 0           & 1           & 1           & 0           & 1           & 1           & 0           & 0          \\
    \end{tabular}
    }
    \label{tab:switchingSecisions}%
    }
\end{table}

\subsubsection{\color{black}{Out-of-sample analysis}} To compare the performance of both solutions provided in Table \ref{tab:switchingSecisions}, we conduct the following out-of-sample analysis. Firstly, we have solved problem \eqref{eq_model4_upper_primal00}--\eqref{eq_model4_upper_primal34} forcing each of the two obtained switching decisions (without considering the last term in the objective function). Given the obtained power flows, we calculated the probability of failure for each line given switching decisions. Then, we generated 2000 scenarios of failure following a Bernoulli trial for the line states (1 in service; 0 failure) with the computed probabilities. Under these generated scenarios, we have evaluated the performances of the two solutions. For this out-of-sample analysis, the average loss of load (\% of total demand) for the solutions {\it without DDU} and {\it with DDU} are 44.15\% and 0.53\%, respectively. In addition, the {\color{black}CVaR$_{95\%}$} of loss of load (\% of total demand) for the solutions {\it without DDU} and {\it with DDU} are 57.17\% and 6.91\%, respectively. Moreover, according to Fig. \ref{Fig:inverseCumulativeDist}, the solution {\it with DDU} has 85.25\% probability to incur in null loss of load and 96.85\% probability to resulting in up to 2\% of loss of load, whereas the solution {\it without DDU} has 98.00\% probability to incur a loss of load and more than 90\% probability to result in more than 30\% of loss of load. Therefore, our proposed model can properly recognize the appropriate switching actions that are needed to significantly decrease the risk of loss of load within a decision-dependent uncertainty framework.

\subsection{138-bus system}

\begin{figure*}[!tb]
    \centering
     \includegraphics[width=0.85\textwidth,height=0.5\textheight,keepaspectratio]{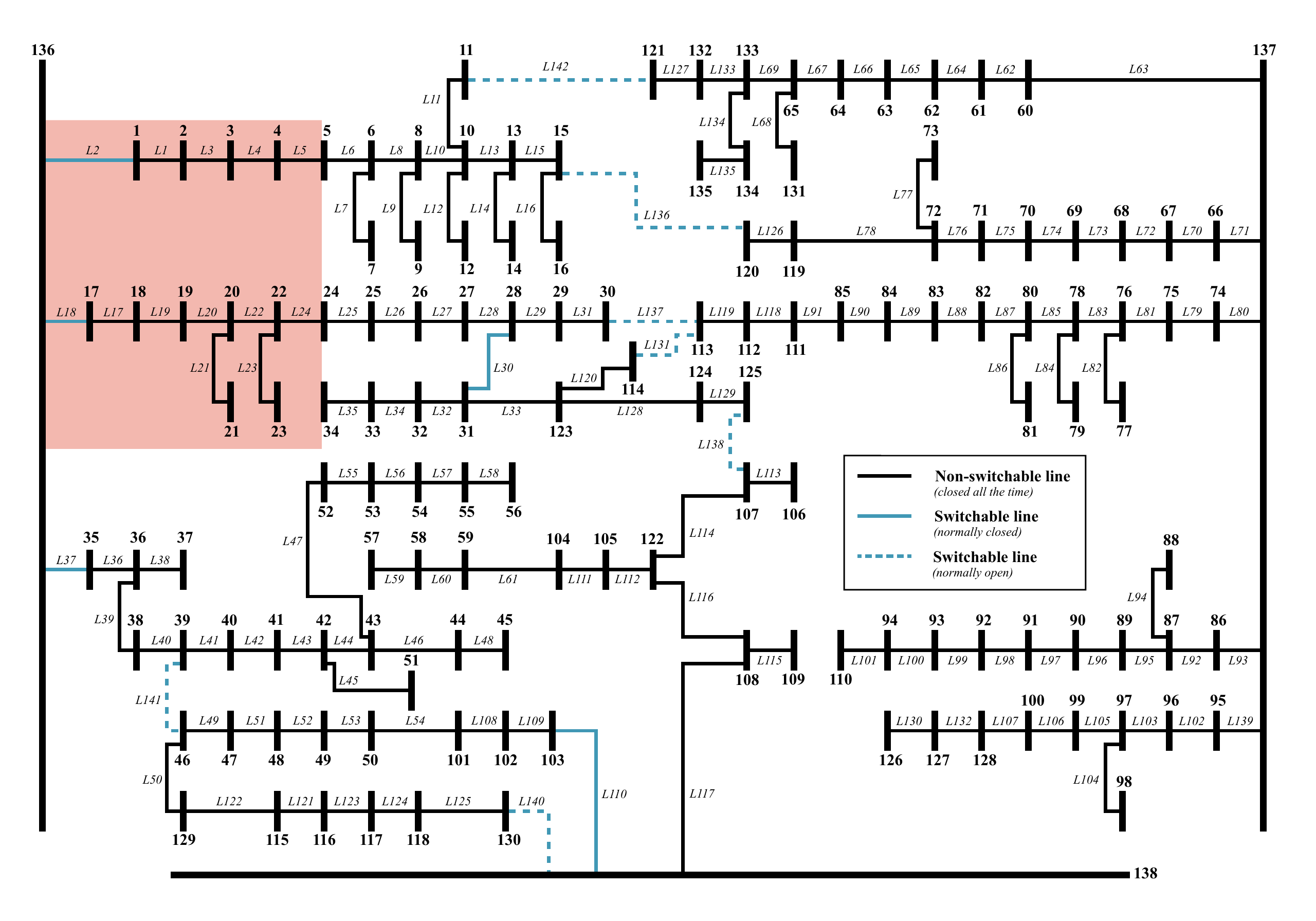}
    \caption{138-Bus distribution system.}
    \label{Fig:138BusFigure}
\end{figure*}

We have also studied the benefits and effectiveness of the proposed methodology in the larger and more complex 138-bus distribution system (Fig. \ref{Fig:138BusFigure}), based on the data provided in \cite{Munoz2016}. In this system, there are 3 substations, 138 buses, and 142 lines, from which 12 are switchable. The total demand of the system is 56,900 kW, the energy price is 0.2 \$/kWh, the deficit cost is 2 \$/kWh, and each switching action costs \$200. To enforce radiality, we used a DFS (depth-first search) algorithm to identify 12 rules that avoid the simultaneous activation of line segments and result in the formation of cycles within the network. All branches have a nominal rate of failure equal to 0.15 per year and, analogously to Subsection \ref{Sec:54bus}, this rate of failure translates into a 0.0411\% of failure probability for each line in the next 24 hours ($\gamma$) using the exponential probability distribution. Furthermore, the northwest part of the system is more likely to initiate a wildfire. This area of the grid includes lines L1--L5 and L17--L24.

\begin{table*}[htbp]
  \centering
  \scriptsize
  \newcolumntype{R}{>{\raggedleft\arraybackslash}m{3.70em}}
  \caption{Main results - 138-bus system}
  \renewcommand{\arraystretch}{1.30}
  \begin{tabular}{>{\raggedleft\arraybackslash}m{16.5em}|R|RRRRRRRR|}
        \multicolumn{1}{r|}{} & \multicolumn{1}{c|}{\textbf{W/out DDU}} & \multicolumn{8}{c|}{\textbf{With DDU}} \\ 
        \hline
            \textbf{Maximum failure probability} & {-} & \textbf{1.0\%} & \textbf{1.1\%} & \textbf{1.8\%} & \textbf{1.9\%} & \textbf{3\%} & \textbf{4\%} & \textbf{50\%} & \textbf{90\%} \\
        \hline
            \textbf{Switching actions (line index)} & \multicolumn{1}{r|}{-} & \multicolumn{1}{r}{-} & \scriptsize{30; 138} & \scriptsize{30; 138} & \scriptsize{2; 30; 136; 138} & \scriptsize{2; 30; 136; 138} & \scriptsize{2; 18; 30; 136; 137; 138} & \scriptsize{2; 18; 30; 136; 137; 138} & \scriptsize{2; 18; 30; 136; 137; 138} \\
        \hline
            \textbf{Objective function value (\$)} & 23,110 & 24,118 & 24,199 & 24,523 & 24,555 & 24,745 & 24,872 & 25,582 & 26,200 \\
            \footnotesize{Energy} & 11,380 & 11,380 & 11,380 & 11,380 & 11,380 & 11,380 & 11,380 & 11,380 & 11,380 \\
            \footnotesize{Switching} & 0 & 0 & 400 & 400 & 800 & 800 & 1,200 & 1,200 & 1,200 \\
            \footnotesize{Deficit} & 0 & 0 & 0 & 0 & 0 & 0 & 0 & 0 & 0 \\
            \footnotesize{Worst case expected value of post-contingency operation cost} & 11,730 & 12,738 & 12,419 & 12,743 & 12,375 & 12,565 & 12,292 & 13,002 & 13,620 \\
        \hline
            \textbf{Computing time (min)} & 2.47 & 7.26 & 8.90 & 14.87 & 15.35 & 15.14 & 18.71 & 16.39 & 19.42 \\
    \end{tabular}%
  \label{tab:main_results_138}%
\end{table*}%

In this numerical experiment, we conducted a sensitivity analysis of the impact of the $\beta$ parameter in the solution by running the model {\it with DDU} 27 times, considering different values for $\beta$ in the mentioned area. The range of the chosen values was defined considering the maximum failure probability ($\gamma + \beta \overline{F}$). This probability indicates how likely a line failure is to happen if the power flow in the feeder is at its maximum capacity. Given that, we chose $\beta$ values for the lines in the wildfire area considering $\beta \times \overline{F}$ to range from 1\% to 2\% by 0.1\%, from 2\% to 10\% by 1\%, and from 10\% to 90\% by 10\%. All the lines from outside the wildfire-prone area were assumed to have a $\beta \times \overline{F}$ as 0.1\% in all cases. The input data can also be downloaded from \cite{CaseStudyRepository_DDU_wildfire2023}. 

The main results are depicted in Table \ref{tab:main_results_138}, where, for expository purposes, only the results for 8 cases are shown. These cases of maximum failure probability are important as they resulted in an operation change in terms of switching actions, for example, the cases between 1.1\% and 1.8\% resulted in the same switching decision, and so on. Besides that, values for the model {\it with DDU} refer to running the {\it with DDU and warm up} setup, since the warm-up helps in decreasing the computational burden to handle the decision-dependent model. 

As depicted in Table \ref{tab:main_results_138}, as the value of $\beta$ increases, the line risk of failure also increases, thus the solution is to change the grid by switching some critical lines. By changing the grid topology, the model decreases the power flow in the critical lines (inside the wildfire-prone area), decreasing the risk of failure associated with $\boldsymbol{\beta}$. Moreover, as we increase the level of $\boldsymbol{\beta}$, the worst-case expected value of post-contingency operation cost increases until it is worth performing switching actions. For instance, in the cases where only 4 lines are switched, between 1.9\% and 3\%, the worst-case expected value increases up to \$12,565. At 4\%, similarly, it is economically viable to afford further two switching actions and have a lower worst-case expected value. In general, the solution time of each case also increases with the $\boldsymbol{\beta}$ levels, reaching a maximum elapsed time of roughly 20 minutes.

\subsubsection{Out-of-sample analysis} 
\begin{figure}[!tb]
    \centering
    \includegraphics[width=.45 \textwidth, height=0.5 \textheight, keepaspectratio]{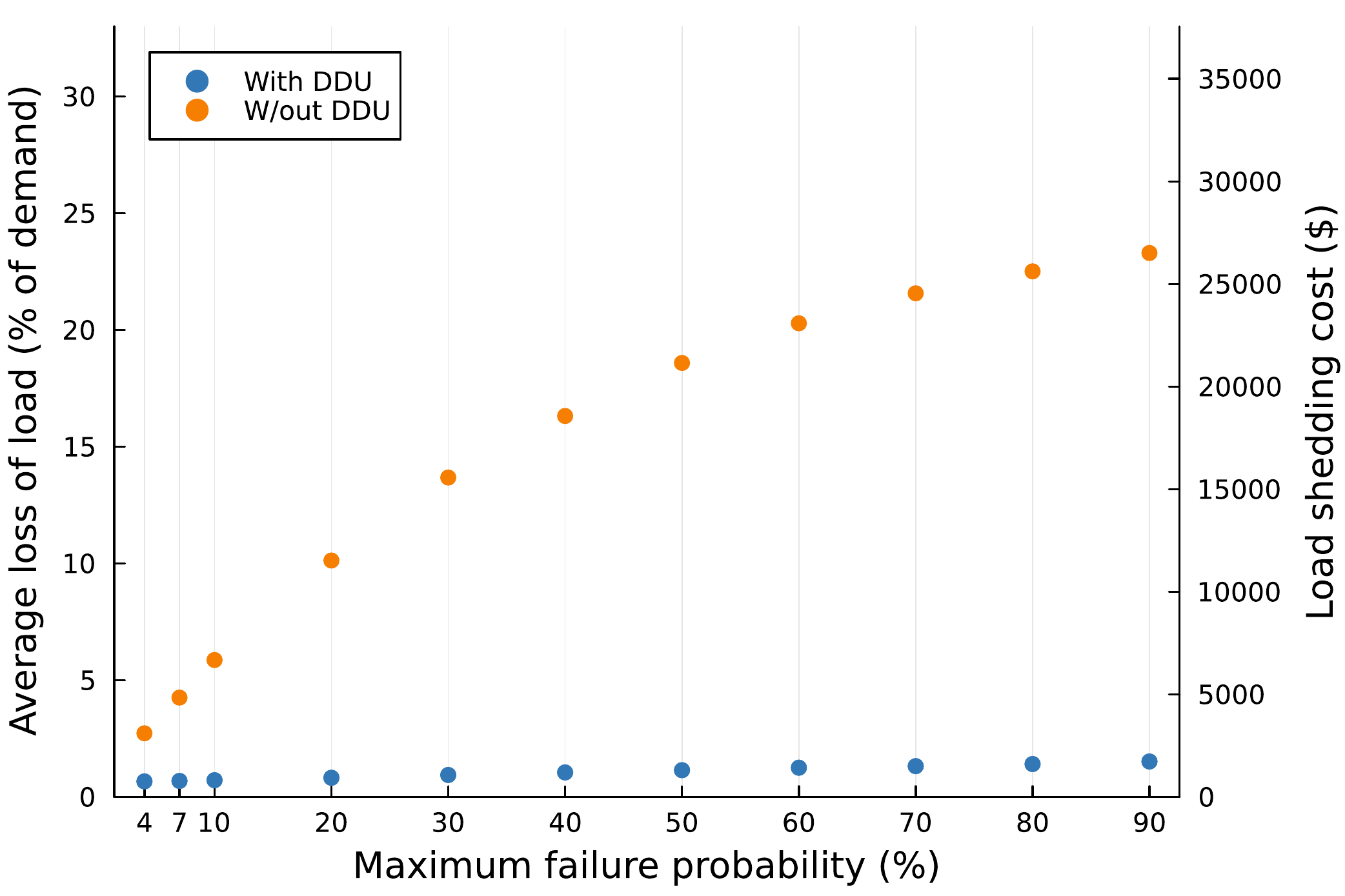}
    \caption{Average loss of load (\% total demand) in the out-of-sample analysis for the solution setup \textit{with DDU} and \textit{without DDU} for different levels of maximum failure probability.}
    \label{Fig:138AverageLoadLoss}
\end{figure}

\begin{figure}[!tb]
    \centering
    \includegraphics[width=.45 \textwidth, height=0.5 \textheight, keepaspectratio]{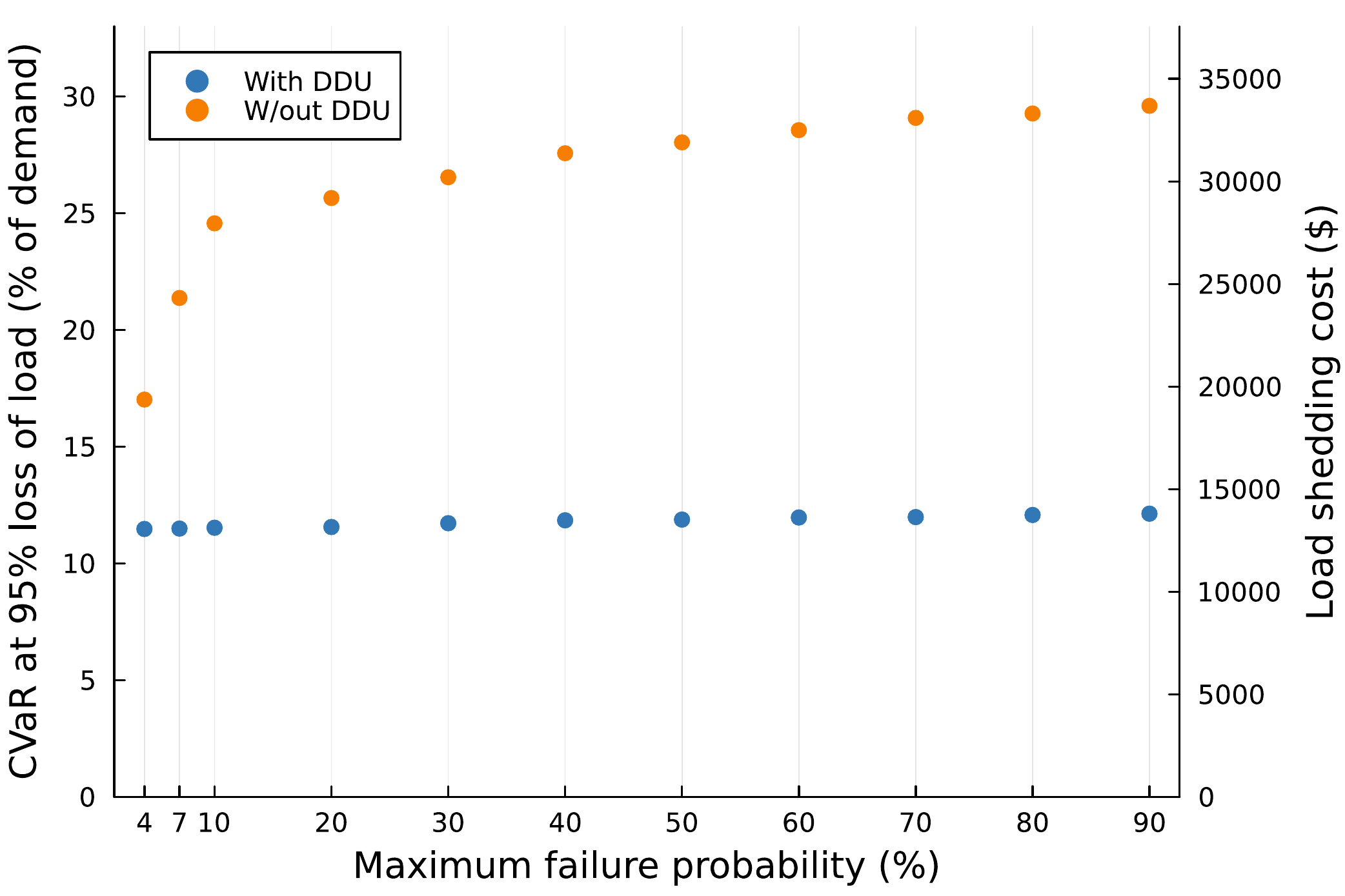}
    \caption{$\text{CVaR}_{95\%}$ loss of load (\% total demand) in the out-of-sample analysis for the solution setup \textit{with DDU} and \textit{without DDU} for different levels of maximum failure probability.}
    \label{Fig:138CVaRLoadLoss}
\end{figure}

Finally, we also performed an out-of-sample analysis using the same procedure presented for the 54-bus system. In the 138-bus system, we consider the results of using each parameter $\boldsymbol{\beta}$. Firstly, in Fig. \ref{Fig:138AverageLoadLoss}, we showcase the average load shedding for the out-of-sample analysis for each level of maximum failure probability ($\gamma + \beta \overline{F}$). Note that, as the environmental conditions for a wildfire worsen, the impact of the average loss of load when disregarding the DDU increases significantly. For instance, for a maximum failure probability of 90\%, the average loss of load would be roughly 22\% of total demand if no actions were considered ({\it without DDU}), while it would be roughly 1\% if the actions suggested by the DDU model were implemented. Furthermore, Fig. \ref{Fig:138CVaRLoadLoss} depicts a similar analysis, but highlighting the associated $\text{CVaR}_{95\%}$ level. Note that, for the setup {\it without DDU} , a value of roughly 30\% in loss of load (in \% of total demand) can be observed in the most critical scenarios. On the other hand, nevertheless, the system topology prescribed by the {\it with DDU} setup significantly mitigates the load shedding occurrence and, consequently, the system operation cost. For instance, consider the maximum failure probability of 90\% once more. The average cost in the {\it without DDU} setup is \$26,512 (Fig. \ref{Fig:138AverageLoadLoss}), which is higher than the expected cost in the 5\% worst-valued scenarios ($\text{CVaR}_{95\%}$), given by \$13,806, when prescribing the network topology based on the {\it with DDU} setup (Fig. \ref{Fig:138CVaRLoadLoss}).

\vspace{-0.3cm}
\section{Conclusion}

This paper proposes a novel methodology to operate distribution systems amidst adverse climate conditions. We acknowledge that the likelihood of a line failure is dependent on its scheduled power flow and aggravated under a wildfire-prone environment. Therefore, in this work, we leverage a Decision-Dependent Uncertainty (DDU) framework to characterize the climate- and power-flow-dependent line availability probability function to devise a wildfire-aware distribution grid operation methodology and prescribe optimal switching actions to decrease the usage level of lines in peril locations, resulting in a more reliable operative condition. Two numerical experiments were conducted to illustrate the effectiveness of the proposed methodology. The results demonstrated that by properly considering DDU, our methodology can keep supplying loads when preventive switching actions are taken. This new configuration leads to a decrease in power flows near the areas where wildfire ignitions are more likely to occur. By doing that, the risk of failure and the risk of load loss are reduced. Although this method can be seen as a better alternative to PSPS, the model can also be adapted to determine where the shut-offs actions should be made.
%
%
\vspace{-0.3cm}
%
\bibliographystyle{IEEEtran}
\bibliography{IEEEabrv,References}

\end{document}